\documentclass[final]{siamltex}

\title{An efficient step size selection for ODE codes}
\author{Lars Petter Endresen\thanks{Institutt for fysikk, 
        NTNU, N-7034 Trondheim, Norway}}

\begin{document}
\maketitle

\begin{abstract}
We give an algorithm for efficient step size control in numerical
integration of non--stiff initial value problems, based on a formula
tailormade to methods where the numerical solution is compared with a
solution of lower order.
\end{abstract}


\section{Introduction}
We discuss here step size control in numerical integration of
ordinary differential equations, for the purpose of optimizing
performance in terms of precision and computer time. Thus the aim is
not to estimate or control absolute errors, although upper limits can
be given. Only non--stiff initial value problems are considered.

A new formula for efficient step size control has recently been
proposed \cite{Endresen3}. This formula has the desirable property
that the value it gives for the next step size $h_{n+1}$ is
independent of the present step size $h_{n}$ in the asymptotic limit
$h_{n}\to 0$. For a one--step, variable step size method of order
$(p,p-1)$ that in addition to the primary integration method of order
$p$, also uses a secondary method of order $p-1$, the step size
control algorithm using this formula can be written:

\begin{eqnarray}
\label{eq:form1}
{\rm   } \qquad {h}_{n+1} &=& h_{n}\left(\tau\over|\epsilon|{h}_{n}\right)^{1\over p+1} \\
{\rm L:} \qquad {\lambda}_{1}h_{n} &\leq& h_{n+1}\leq {\lambda}_{2} h_{n}\\
{\rm A:} \qquad |\epsilon| h_{n} &<& \sigma \tau 
\end{eqnarray}
 
\noindent
Here $\epsilon$ is the difference between the two numerical solutions
in one single step of length $h_{n}$ and $\tau$ is a constant tolerance.
We limit the increase and decrease of the step size by the limitation
criterion L, and reject steps that do not satisfy the acceptance
criterion A. $\sigma$, ${\lambda}_{1}$ and ${\lambda}_{2}$ are
parameters. We show that the above algorithm is more efficient than
the standard algorithm known from the literature
\cite{HairerNorsettWanner}:

\begin{eqnarray}
\label{eq:form2}
{\rm   } \qquad h_{n+1} &=& h_{n}\left(\tau\over|\epsilon|\right)^{1\over p+1}\\
{\rm L:} \qquad {\lambda}_{1}h_{n} &\leq& h_{n+1}\leq {\lambda}_{2} h_{n}\\
{\rm A:} \qquad \quad |\epsilon| &<& \sigma \tau 
\end{eqnarray}

\noindent
One way to measure efficiency for one--step methods is to determine
the number of function evaluations for a given global error. We use
the software package DETEST \cite{Hull}, applying the embedded
Runge--Kutta method of Dormand and Prince \cite{DormandPrince}, of
order $(p,q)=(5,4)$. First we find the values of $\sigma$,
${\lambda}_{1}$ and ${\lambda}_{2}$ that give the least number of
function evaluations for a global error of $10^{-4}$, $10^{-5}$,
$10^{-6}$ and $10^{-7}$. This is done simultaneously for 10 of the
problems (group I) in DETEST using a Levenberg--Marquardt method with
a mixed quadratic and cubic line search procedure \cite{Matlab}. We
then compare the two algorithms with this set of parameters on 10
other problems (group II) from DETEST. The computer used was a Cray
Origin 2000.

\section{Numerical Results}
Group I is defined arbitrarily to be the problems A1, A3, A5, B2, B4,
C1, C3, C5, E2 and E4, and group II to be the problems A2, A4, B1, B3,
B5, C2, C4, E1, E3 and E5, in DETEST \cite{Hull}. We omitted the
problems in class D and F, since DETEST only gave results for a
limited set of tolerances in these (the tolerance is varied
automatically by the program). The parameters that resulted in the
least number of function evaluations for group I were:

\begin{equation}
\label{eq:form4}
\begin{array}{llll} 
\sigma = 6.70 & {\lambda}_{1} = 0.67 & {\lambda}_{2} = 5.00 & {\rm (formula \; \ref{eq:form1} )}, 
\end{array} 
\end{equation} 
\begin{equation}
\label{eq:form5}
\begin{array}{llll} 
\sigma = 5.50 & {\lambda}_{1} = 0.26 & {\lambda}_{2} = 4.00 & {\rm (formula \; \ref{eq:form2})}. 
\end{array} 
\end{equation} 

\noindent
The new formula gave a lower number of function evaluations in
$50.9$\% of the cases in problem group I (for which $\sigma$,
${\lambda}_{1}$ and ${\lambda}_{2}$ were optimized), and in $64.4$\%
of the cases in problem group II. Only calculations giving global
errors in the range $(10^{-3},10^{-8})$ were counted, since DETEST did
not give results outside this range for some problems. In Table
\ref{table1} we have displayed the mean ratio of the number of
function evaluations of the new formula to the number of function
evaluations of the standard formula. The standard formula was tested
with the optimized set of parameters (equation \ref{eq:form4}) and
the recommended set of parameters:

\begin{equation}
\label{eq:form6}
\begin{array}{lll} 
\sigma = 1.20 & {\lambda}_{1} = 0.50 & {\lambda}_{2} = 2.00, 
\end{array} 
\end{equation} 

\begin{table}[h]
\caption{The mean ratio of the number of function evaluations of the
new formula to the number of function evaluations of the standard
formula.}
\begin{center} \footnotesize
\begin{tabular}{|c|c|c|c|c|} \hline 
\multicolumn{1}{|c|}{Expected} &
\multicolumn{4}{|c|}{Mean ratio of function evaluations} \\
\multicolumn{1}{|c|}{global} &
\multicolumn{2}{|c|}{with equation \ref{eq:form4} and \ref{eq:form5} } &
\multicolumn{2}{|c|}{with equation \ref{eq:form4} and \ref{eq:form6} } \\
error & \hspace{0.6cm} Group I \hspace{0.6cm} & \hspace{0.6cm} Group II \hspace{0.6cm} & \hspace{0.6cm} Group I \hspace{0.6cm} & \hspace{0.6cm} Group II \hspace{0.6cm} \\ \hline 
$10^{-3}  $ & $  0.9575  $ & $ 1.0103 $ & $  0.8794  $ & $ 0.8747 $ \\ 
$10^{-4}  $ & $  1.0152  $ & $ 1.0448 $ & $  0.9222  $ & $ 0.9011 $ \\  
$10^{-5}  $ & $  1.0166  $ & $ 0.9696 $ & $  0.9281  $ & $ 0.9266 $ \\  
$10^{-6}  $ & $  1.0014  $ & $ 0.9615 $ & $  0.9370  $ & $ 0.8797 $ \\ 
$10^{-7}  $ & $  1.0000  $ & $ 0.9518 $ & $  0.9397  $ & $ 0.8705 $ \\  
$10^{-8}  $ & $  0.9947  $ & $ 0.9545 $ & $  0.9434  $ & $ 0.8814 $ \\ \hline
\end{tabular}
\end{center} 
\label{table1} 
\end{table}

\section{Discussion}
A new formula for step size selection in numerical integration of
non--stiff initial value problems has been tested on 20 initial value
problems. It is found that this formula on the average is more
efficient than the standard step size selection formula. I thank the
authors of DETEST.

\end{document}